\documentclass[10pt]{article}

\usepackage[T1]{fontenc}
\usepackage[english]{babel}
\usepackage{graphicx}

\usepackage{geometry}
\usepackage{movie15}

\usepackage{animate}
\usepackage{ifdraft}
\usepackage{amsmath}
\usepackage{amsfonts}
\usepackage{amsthm}
\usepackage{amssymb}
\usepackage{amscd}
\usepackage{epsfig}
\usepackage{verbatim}
\usepackage{fancybox}
\usepackage{moreverb}
\usepackage{psfrag}
\usepackage{cellspace}
\usepackage{hyperref}
\usepackage{latexsym}

\usepackage{color}
\definecolor{vert}{rgb}{0.,0.5,0.}
\definecolor{rouge}{rgb}{0.8,0.,0.}
\definecolor{violet}{rgb}{0.5,0.,0.4}
\definecolor{bleu}{rgb}{0.,0.,0.5}
\definecolor{orange}{rgb}{0.8,0.4,0.}
\definecolor{light-blue}{rgb}{0.5,0.5,0.7}
\definecolor{light-red}{rgb}{0.8,0.2,0.2}
\definecolor{noir}{rgb}{0.,0.,0.}
\newcommand{\N}{{\mathbb N}}

\newcommand{\Ea}{{\mathcal E}}

\newcommand{\Pa}{{\mathcal P}}
\newcommand{\Wa}{{\mathcal W}}
\newcommand{\WaFlow}{{\mathcal F}_{\Wa}}
\newcommand{\WaRate}{{\omega}}
\newcommand{\coefWaRate}{{\beta_0}}
\newcommand{\coefWaToEnt}{{\beta_1}}
\newcommand{\coefWaToEff}{{\beta_2}}
\newcommand{\coefWaToPop}{{\beta_3}}
\newcommand{\coefDiffWa}{{\nu}}
\newcommand{\EnToWa}{{\kappa}}
\newcommand{\VectNorm}{{\mu}}

\newcommand{\disk}{{\tt D}}

\newcommand{\imf}{{\mathfrak i}}
\newcommand{\HomeAttraction}{{\mathcal A_{\Pa}}}
\newcommand{\EnterpriseAttraction}{{\mathcal A_{\Ea}}}
\newcommand{\EnterpriseWeight}{{\alpha}}
\newcommand{\MorningDisplacementTime}{{\mathfrak{t}_{m}}}
\newcommand{\EveningDisplacementTime}{{\mathfrak{t}_{e}}}

\newcommand{\fracp}[2]{\frac{\partial #1}{\partial #2}}

\newcommand{\ds}{\displaystyle}

\begin{document}

\title{A PDE-like Toy-Model of Territory Working
\footnote{This work is supported by \href{http://www.mgdis.fr/}{MGDIS} and \href{http://www.agence-maths-entreprises.fr/a/?q=fr}{AMIES}}}
\author{ {\bf Emmanuel Frénod\footnote{
Univ. Bretagne - Sud,  UMR 6205, \href{http://www.lmba-math.fr/}{LMBA}, F-56000 Vannes, France.
\href{mailto:emmanuel.frenod@univ-ubs.fr}{emmanuel.frenod@univ-ubs.fr}.
\href{http://web.univ-ubs.fr/lmam/frenod/index.html}{\url{http://web.univ-ubs.fr/lmam/frenod/index.html}}.
}}}
\date{\empty}

\maketitle

{ \scriptsize
{\bf Abstract - }
This note introduces a PDE-like Toy-Model that embeds several aspects of Territory Working. The long term goal is to build a software tool that behaves like a Territory, and in particular that incorporates its multi-scale-in-time-and-space nature, in order to make simulations, to explore scenarios, to foresee policy impacts, to help to make a decision when facing a change in the environment or in cultural behavior, {\it etc.}. The aim of this note is prove that the concept of building a model that couple Systemic Approach and PDE tools to achieve the evoked long term goal is possible.
}
\section{Introduction}
This note introduces a PDE-like Toy-Model that embeds several aspects of Territory Working. The long term goal is
to build a software tool that behaves like a Territory, and in particular that incorporates its multi-scale-in-time-and-space nature, 
in order to make simulations, to explore scenarios, to foresee policy impacts,
to help to make a decision when facing a change in the environment or in cultural behavior, {\it etc.}.
The aim of this note is prove that the concept of building a model that couple Systemic Approach and PDE tools to achieve the
evoked long term goal is possible.

\begin{figure}[ht]
\begin{center}
\includegraphics[height=6cm, width=10cm]{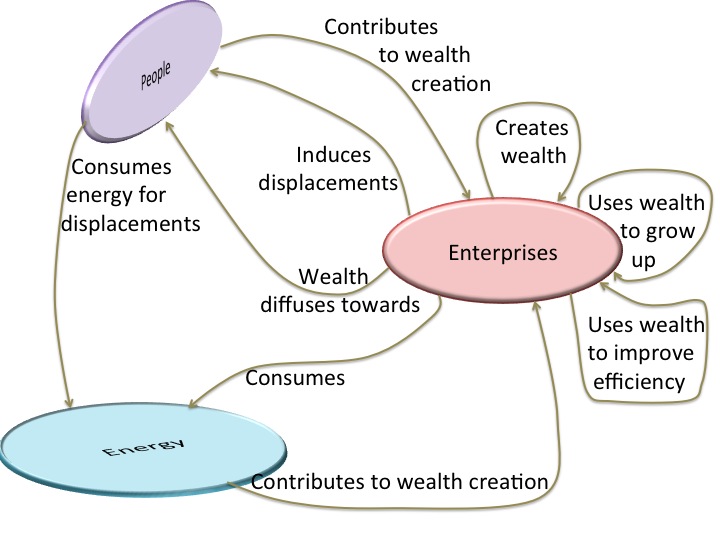}
\caption{Systemic Approach of Territory Working} \label{201212160858} % \ref{201212160858}
\end{center}
\end{figure}
\section{Brief Systemic Approach of Territory Working}
By Territory, we mean a Country or a Town or a set of Towns and Countries having a coherence
and being administrated by a common local government.\\

The Woking of a Territory results from the interactions between several compartments. A representation of
this Systemic Approach is given in Figure \ref{201212160858}. Three compartments are considered. The first one
("People") concerns the population, the second one ("Energy") concerns the energy questions  and the last one
("Enterprises") concerns the world of enterprises . Then, the compartments influence each other or themselves. This fact is 
symbolized by the arrows.  For instance, the arrow which is the more at the top of the figure translates that when people
are at work, they contribute to the production of wealth. In return, the arrow just below means that displacements
are induced by enterprise locations and the arrow in the left translates that an energy consumption is induced by
those displacements.
The three arrows that point on the "Economy" compartment from itself symbolize that enterprises create wealth and
that this created wealth is used by the enterprises for growing up and for efficiency improvement.
The two last arrows, in the bottom, express that enterprises consume energy to create wealth. 
\begin{figure}[ht]
\begin{center}
\includegraphics[height=6cm]{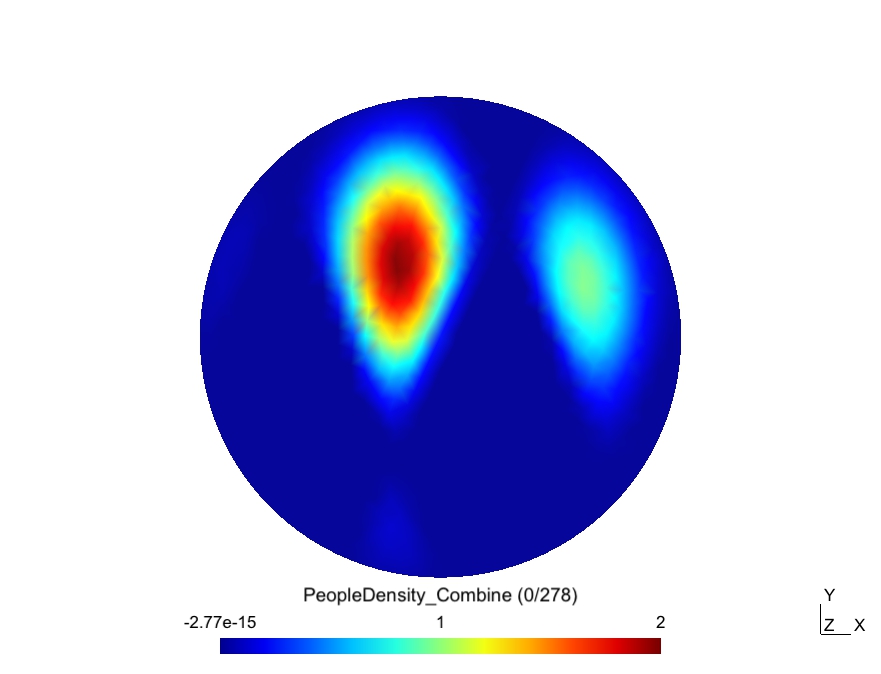}
\hspace{-10mm}
\includegraphics[height=6cm]{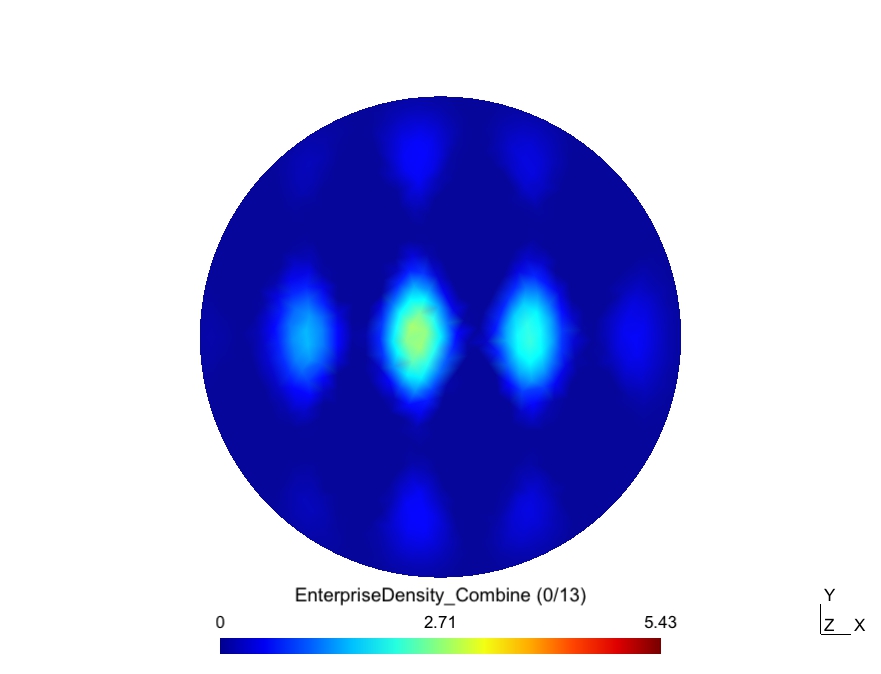}
\caption{Initial Population Density $\Pa_0$ (left) and Job-station Density  $\Ea_0$ (right)}
\label{FigDensInit} % \ref{FigDensInit}
\end{center}
\end{figure}
\begin{figure}[ht]
\begin{center}
\includegraphics[height=6cm]{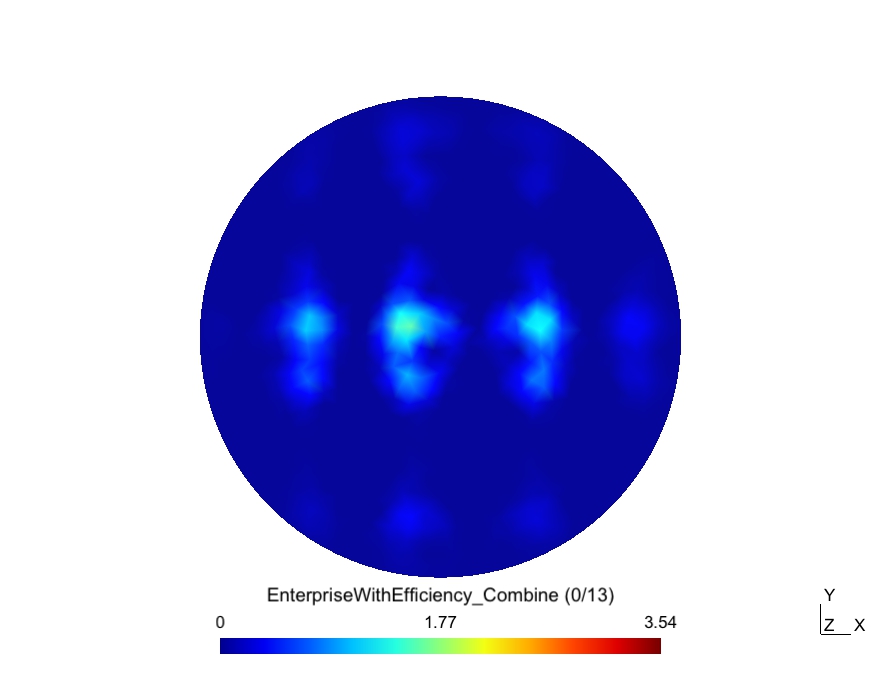}
\caption{Job-station Density $\Ea_0$  times Efficiency Indicator $\imf_0$}
\label{NONFigDensInit} % \ref{NONFigDensInit}
\end{center}
\end{figure}
\begin{figure}[ht]
\begin{center}
\includegraphics[height=6cm]{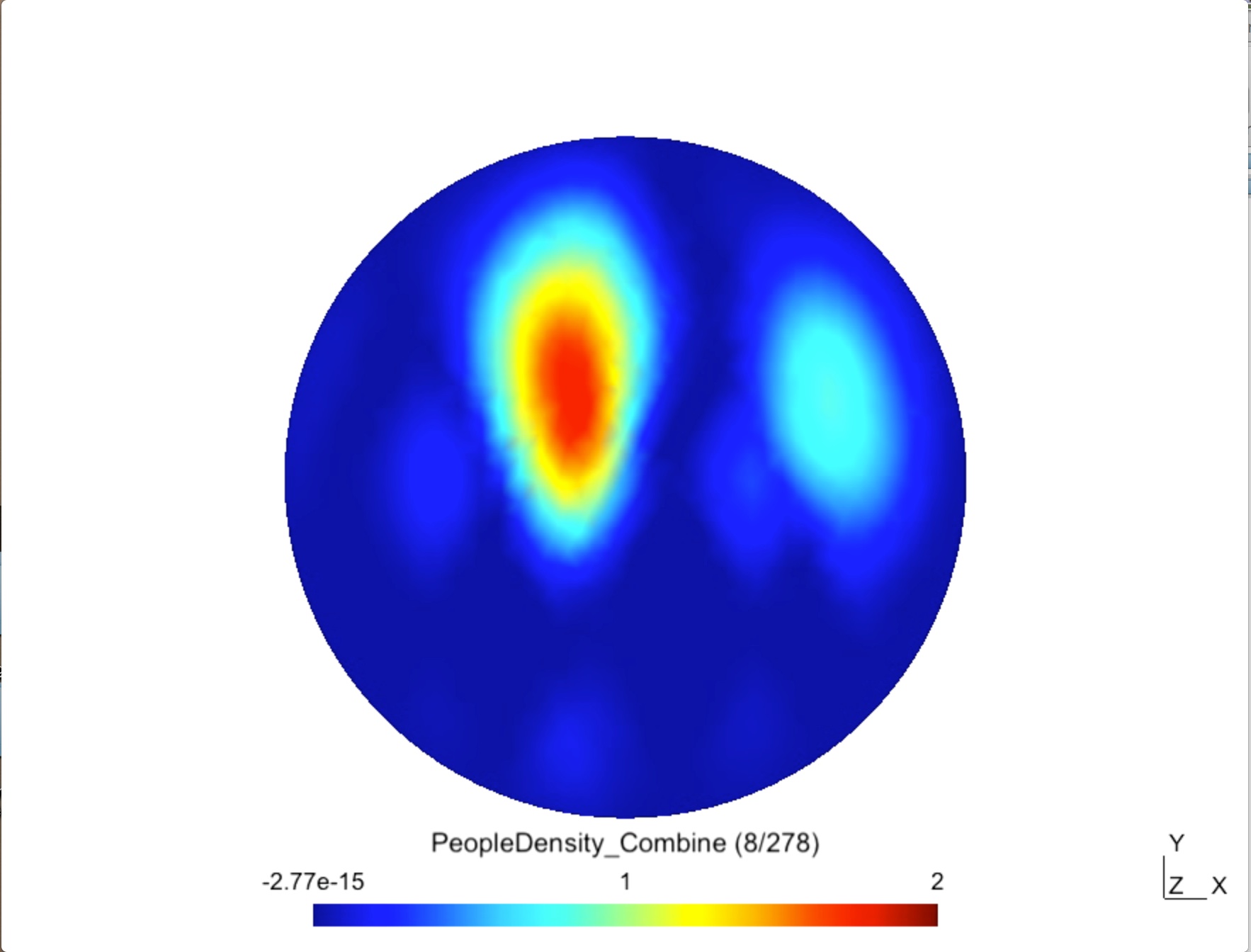}
\caption{Displacements from home to work and back (Animation available on YouTube: \href{http://youtu.be/LkPlVT-a4pg}{\url{http://youtu.be/LkPlVT-a4pg}}.)} 
\label{DailyDispl}  % \ref{DailyDispl}
\end{center}
\end{figure}
\section{The Toy-Model}
In this Toy-Model, all quantities are dimensionless.
The simulations are done using  \href{http://www.freefem.org/ff++/}{Freefem++}. All the simulations are done over the same time interval which is run by variable $t\in[0,T)$
for a real number $T>0$.
Yet as the characteristic times of the various phenomena are not the same, they are given with various steps.   \\

The geographic model of the Toy-Model is a disk $\disk$,  provided with coordinates $x$, with a boundary $\partial \disk$. The vector of norm 1, orthogonal to  $\partial \disk$
and pointing outside  $\disk$ is denoted $\VectNorm$.

On this disk, several densities, which depend on time, are defined: 
\begin{itemize}
\item $\Pa=\Pa(t,x)$ is the Population Density.
\item $\Ea=\Ea(t,x)$ is the Job-station Density.
\item $\Wa=\Wa(t,x)$ is the Wealth-that-goes-to-people Density.
\end{itemize}
\begin{figure}[ht]
\begin{center}
\includegraphics[height=6cm]{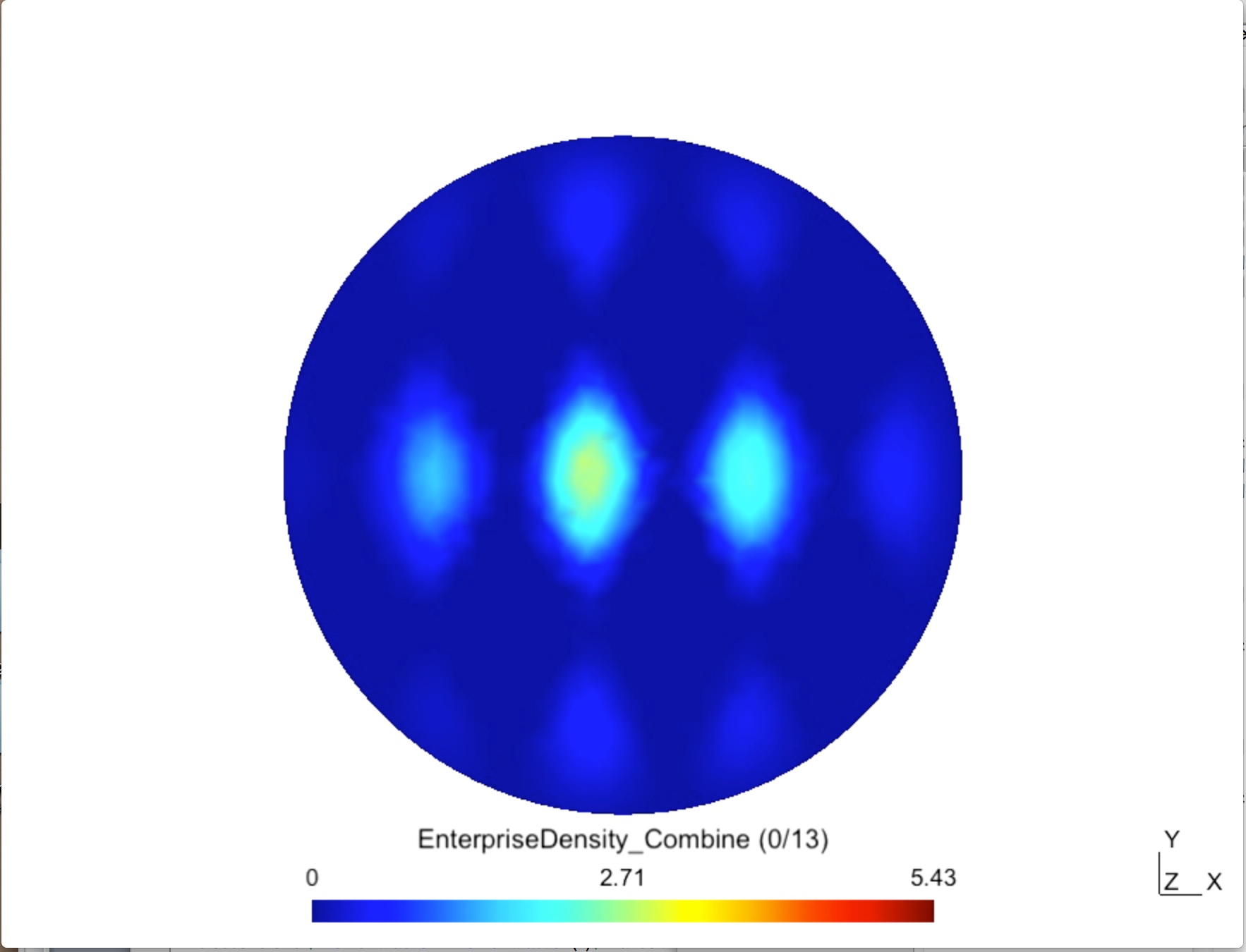}
\caption{Evolution of Job-station  Density $\Ea$. (Animation available on YouTube: \href{http://youtu.be/PorgK4q5Ds0}{\url{http://youtu.be/PorgK4q5Ds0}}.)}
\label{201212121404} % \ref{201212121404}
\end{center}
\end{figure}
For those densities, initial values are defined: $\Pa_0=\Pa_0(x)$, which is drawn on the left of Figure \ref{FigDensInit}, and which is the Population Density when everybody is at home at night; $\Ea_0=\Ea_0(x)$ which is drawn on the right of Figure \ref{FigDensInit}; and, $\Wa_0=\Wa_0(x)=1$.

An Efficiency Indicator $\imf=\imf(t,x)$ is also defined on the disk. It ranges in $[0,1]$ and measures the efficiency of the enterprises at every point of the disk; it is initialized at $\imf_0$. The product of the initial Job-station Density $\Ea_0$ by the initial Efficiency Indicator $\imf_0=\imf_0(x)$  is given in Figure \ref{NONFigDensInit}. \\

From  initial Population Density $\Pa_0$ and Job-station Density $\Ea$ two attractors are built. 
The first one $\EnterpriseAttraction=\EnterpriseAttraction(t,x) $ attracts the Population Density towards Job-stations.
It is defined as
\begin{gather}
  \EnterpriseAttraction(t,x) =
  \frac{ \ds \int_\disk \Pa_0(x) \, dx}{\ds \int_\disk \Pa_0(x) \, dx + \EnterpriseWeight \int_\disk \Ea(t,x) \, dx} \left(  \Pa_0(x) + \EnterpriseWeight \Ea(t,x) \right),
\end{gather}
where $\EnterpriseWeight$ is a constant larger than 1.
The second one $\HomeAttraction=\HomeAttraction(x)$ makes the population to go back home. It is defined as
\begin{gather}
\HomeAttraction(x) = \Pa_0(x).
\end{gather}
 
With those two attractors, we can simulate the daily motion of the population. 
For this, we firstly adopt the convention that day length is 1 and we define two functions $\MorningDisplacementTime=\MorningDisplacementTime(t)$
and $\EveningDisplacementTime=\EveningDisplacementTime(t)$ depending only on time and periodic of period 1.
$\MorningDisplacementTime(t)>0$ for instants corresponding to the morning displacements, when people go to work and 0 otherwise; 
and, $\EveningDisplacementTime(t)>0$ for instants corresponding to the evening displacements, when people go back home and 0 otherwise.
Secondly we write the equation meaning that in the morning the Population Density is attracted by Job-station Attractor and in the evening by Home Attractor:
\begin{gather}
\fracp{\Pa}{t}(t,x) =\MorningDisplacementTime(t) \left( \EnterpriseAttraction(t,x)-\Pa(t,x)  \right)+  \EveningDisplacementTime(t) \left(\HomeAttraction(x)-\Pa(t,x)\right), ~\forall x\in\disk, \forall t\in(0, +\infty),
\end{gather}
The evolution of ${\Pa}$ is given by the movie in Figure  \ref{DailyDispl}.
In this movie, we can see the alternation of people displacements to their works (in the morning) and to their homes (in the evening).
Because of Job-station distribution most of people converge everyday to work in a small region located near the center of the disk.\\

This motion generates an energy consumption. It is modeled by a time density $\phi=\phi(t)$, which is constant per day and which value at any given day is the  
double integral over the disk of the distances from home locations to Job-station locations weighted by the locations where Population Density $\Pa$ increases and where it decreases.\\
In other words, $\phi$ is a constant over every interval $[n,n+1], n\in\N$ with worth
\begin{gather}
 \phi^n=\int_{\disk}\int_{\disk} |x-y| \,\Pa^n_\text{\rm Incr}(x) \Pa^n_\text{\rm Decr}(y) \, dx dy,
\end{gather}
where $\Pa^n_\text{\rm Incr}$ is the density of increasing population, where population increases, when comparing the Population Density in the morning 
(before people go to work) and in the middle of the day (when workers are at their job-station). It is defined by
\begin{gather}
\Pa^n_\text{\rm Incr}(x) = \max(\Pa(n+\frac 12, x) - \Pa(n, x), 0).
\end{gather}
In a similar way, $\Pa^n_\text{\rm Decr}$ is the density of decreasing population, where population decreases, when comparing the Population Density in the morning and in the middle of the
day and is defined by
\begin{gather}
\Pa^n_\text{\rm Decr}(x) = \max(-\Pa(n+\frac 12, x) + \Pa(n, x), 0).
\end{gather}

The fact that people is at work produces Wealth with a rate which depends on the product of the Population Density $\Pa$ by the Job-station Density  $\Ea$ times the Efficiency Indicator $\imf$. 
To model this, we introduce a wealth production rate density $\WaRate= \WaRate(t,x)$ which is constant per day and which value at any given day, represented by interval $[n,n+1], n\in\N$, is
\begin{gather}
\WaRate^n(x) =  \coefWaRate\,\Pa(n+\frac 12, x)\,  \Ea(n+\frac 12, x)\,\imf(n+\frac 12, x).
\end{gather}
where $\coefWaRate$ is a coefficient small in front of $1$, meaning that the time scale of Wealth variation is large when compared to a day.

This Wealth is split into three parts.  The first one is allocated to enterprises' growth. In other words, it is a source term of the ODE, set in every point of the disk, the Job-station Density $\Ea$ is solution to:
\begin{gather}
\fracp{\Ea}{t} (t,x)= \coefWaToEnt \, \WaRate(t,x), ~\forall x\in\disk, \forall t\in(0, +\infty),
\end{gather}
where $\coefWaToEnt$ is a coefficient belonging to a triplet $(\coefWaToEnt, \coefWaToEff, \coefWaToPop)$ of positive coefficients  such that $\coefWaToEnt + \coefWaToEff +\coefWaToPop =1$.

The second part ($\coefWaToEff \, \WaRate(t,x)$) contributes to improve enterprises' efficiency. It contributes as a source term in the ODE, set in every point of the disk and which is such that its solution $\imf$ always ranges in $[0,1]$:
\begin{gather}
\fracp{\imf}{t} (t,x)= \imf(t,x) (1-\imf(t,x)) \left(\coefWaToEff \, \WaRate(t,x)- \coefWaToEnt \, \WaRate(t,x)\right), ~\forall x\in\disk, \forall t\in(0, +\infty).
\end{gather}
The term $\coefWaToEnt \, \WaRate(t,x)$ in this equation models that when organization grows, if nothing is done, its efficiency decreases.

\begin{figure}[ht]
\begin{center}
\includegraphics[height=6cm]{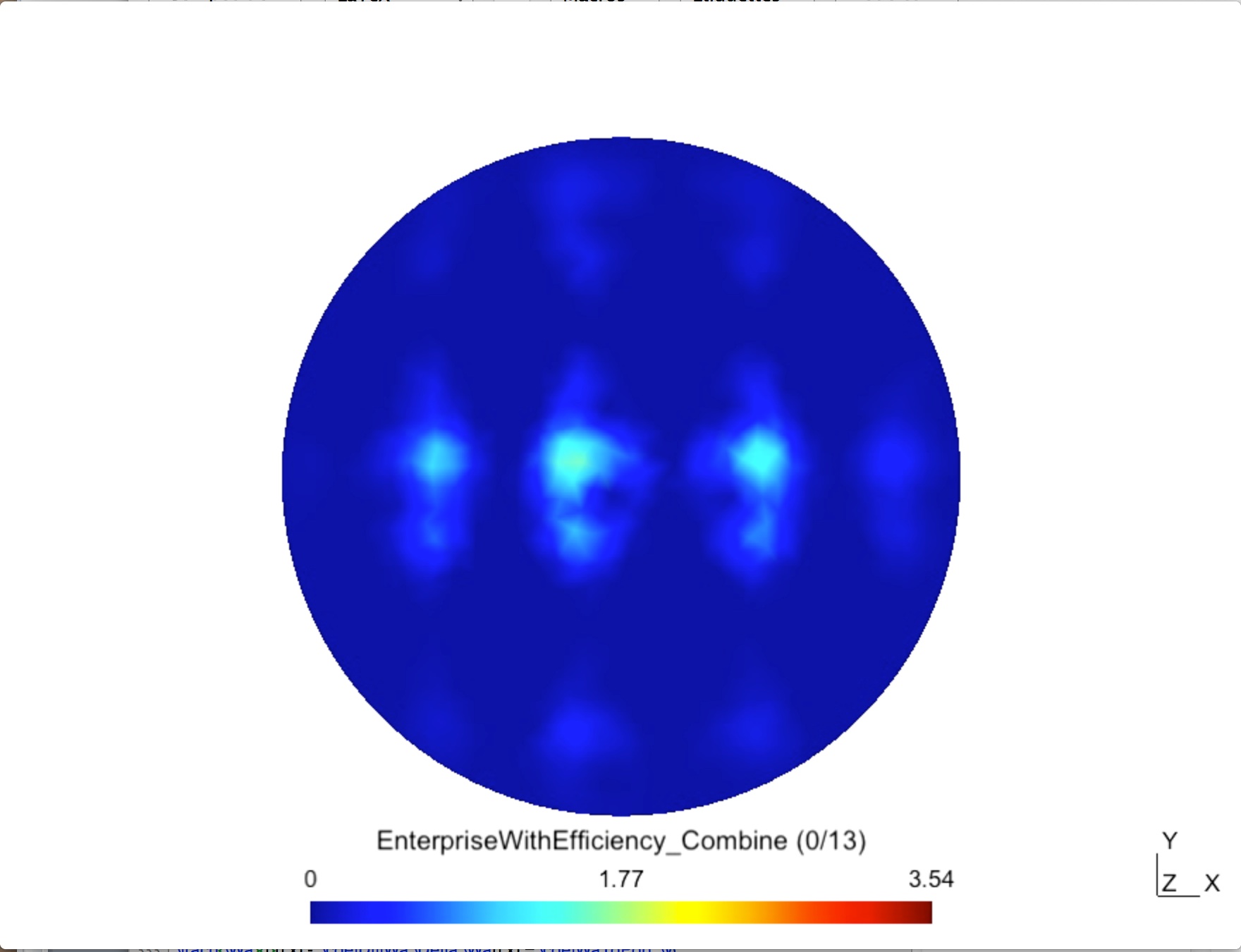}
\caption{Evolution of Job-station  Density $\Ea$ times Efficiency Indicator $\imf$. (Animation available on YouTube: \href{http://youtu.be/B7gjmu7q8IU}{\url{http://youtu.be/B7gjmu7q8IU}}.)}
\label{201212121403} % \ref{201212121403}
\end{center}
\end{figure}

The result of the influence of Wealth on the Job-station Density $\Ea$ and on the product of the Job-station Density  $\Ea$ times the Efficiency Indicator $\imf$ are given in the movies of Figures  \ref{201212121404} and \ref{201212121403}. 
We see in those movies that,  around the small region where most of people converges everyday to work, there is a strong increase of Job-station Density  $\Ea$ and of enterprises' Efficiency  $\imf$.
\begin{figure}[ht]
\begin{center}
\includegraphics[height=6cm]{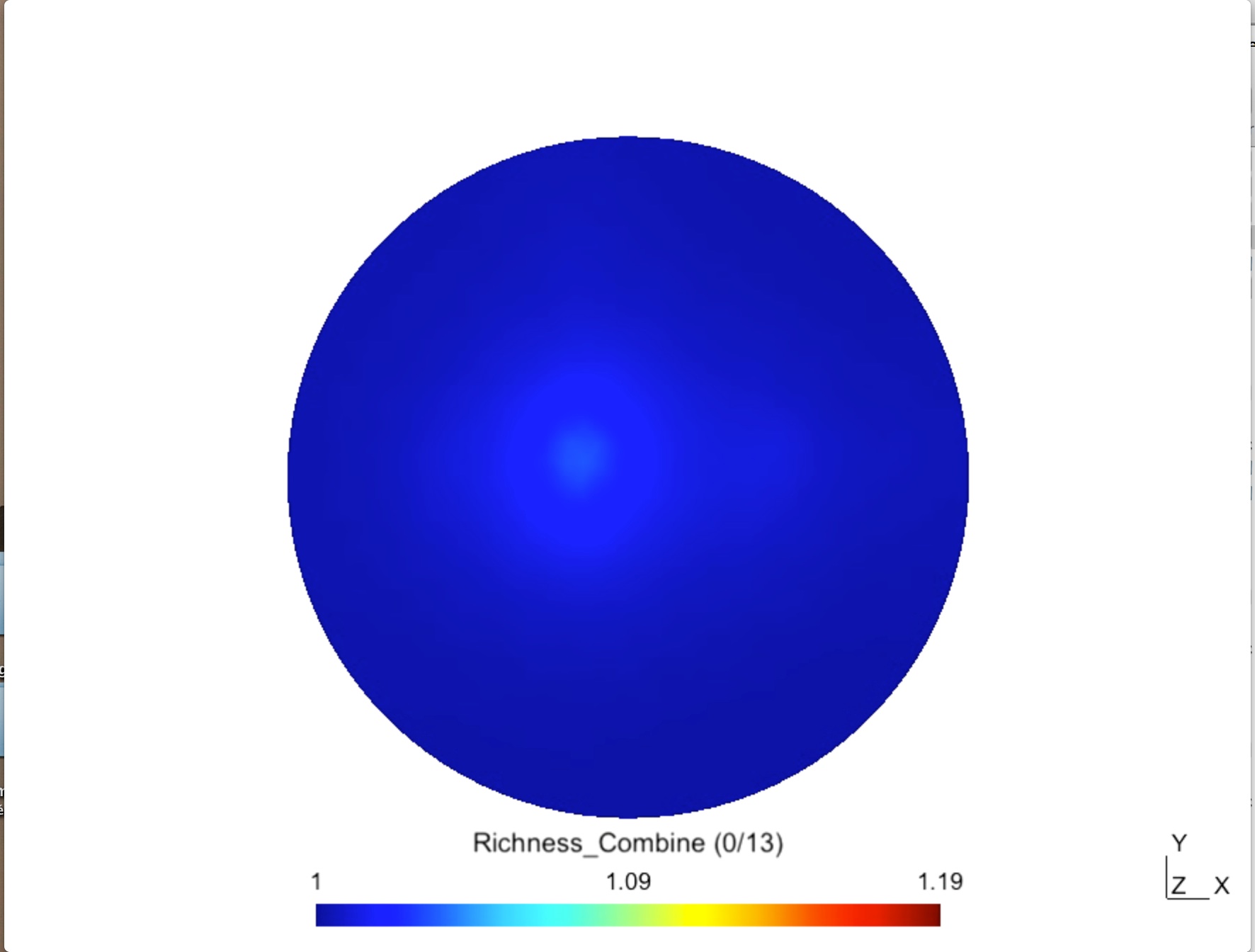}
\caption{Evolution of Wealth Density $\Wa$. (Animation available on YouTube: \href{http://youtu.be/05JNyznPxvw}{\url{http://youtu.be/05JNyznPxvw}}.)}
\label{20121212281612} % \ref{20121212281612}
\end{center}
\end{figure}

The third part of Wealth production will diffuse inside the population. It is then a source term in the diffusion equation $\Wa$ is solution to:
\begin{gather}
\fracp{\Wa}{t}(t,x) -  \coefDiffWa \Delta \Wa(t,x) = \coefWaToPop  \WaRate(t,x) - \EnToWa \phi(t), ~\forall x\in\disk, \forall t\in(0, +\infty),
\\
\fracp{\Wa}{\VectNorm} = \WaFlow, ~ \forall x\in\partial \disk, \forall t\in(0, +\infty),
\end{gather}
where $\coefDiffWa$ is the diffusion coefficient of Wealth within the population on the territory.
In this PDE, the energy consumption time density $\phi$ is also a source term but with an opposite action than that of the Wealth production. This is because $\phi$ induces a consumption of Wealth. In this PDE there is also a boundary condition on the disk border, translating the fact that there is a little wealth that enters the Territory through its border.
The evolution of Wealth-that-goes-to-people Density is given by the movie in Figure  \ref{20121212281612}.
We can see the diffusion of the  Wealth-that-goes-to-people from the small region where most of people converges everyday to work towards the other parts of the Territory.

\section{Interpretation}
Clearly, this Toy-Model does not pretend to give realistic results. Nevertheless it seems that the proposed approach brings a way to couple  several aspects of Territory Working.\\

\begin{figure}[h]
\begin{center}
\includegraphics[height=6cm]{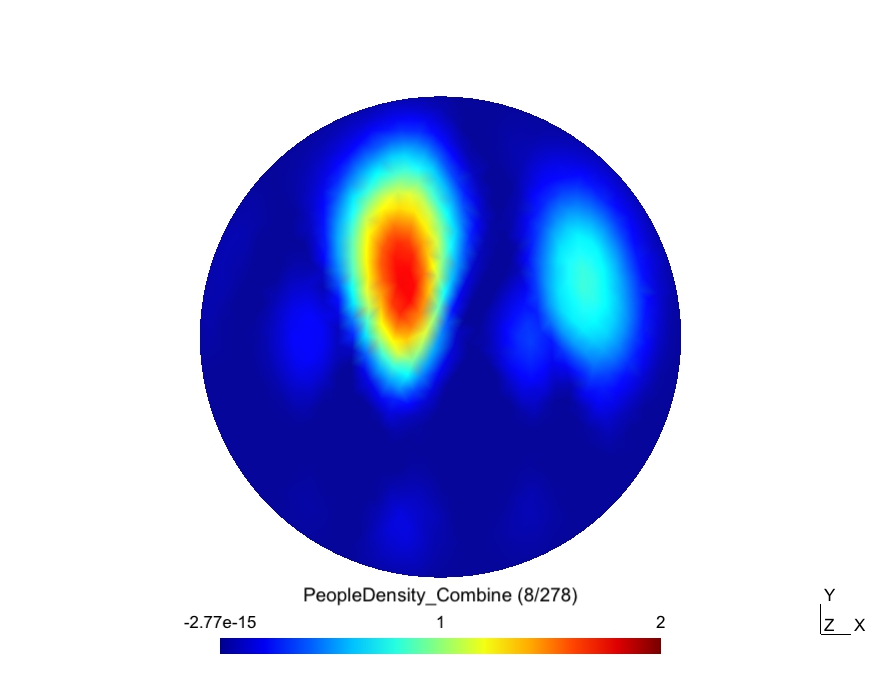}
\hspace{-10mm}
\includegraphics[height=6cm]{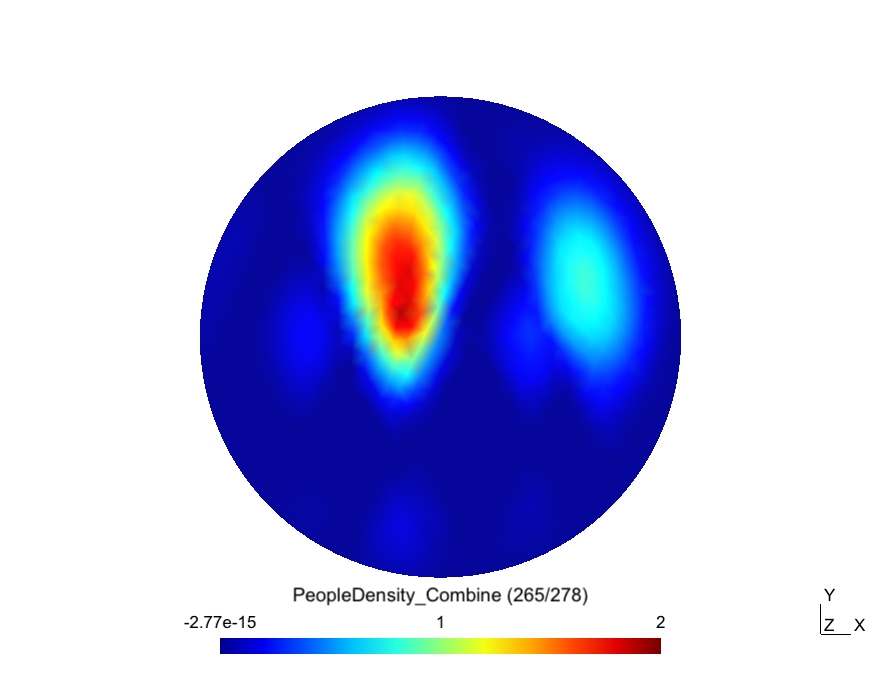}
\caption{Population Density at midday of the first day (left) and of the last day (right)}
\label{20121212281630} % \ref{20121212281630}
\end{center}
\end{figure}

For instance, it can handle non-linearities. This is illustrated by the two pictures of Figure \ref{20121212281630}. The one on the left is the Population Density at midday of the first day and  the one on the right is the Population Density at midday of the last day. 
We can see that, in the small region close to the center of the disk, the density is higher in the right picture than in the left picture. 
This can be explained as follow: Since people is working their, wealth is produced, which induces Job-station number to increase.
As a consequence, day after day, more and more people is coming to work in this small region.\\

The scales (regarding time and space) have no realistic meaning. Nevertheless, we can see the capability of this kind of models to account for a wide variety of scales.
For instance, regarding the space scales, in the Toy-model, there is the size of the Territory (the disk) and there are the characteristic sizes of variation of Population-at-home Density, of Job-station Density and Enterprise Efficiency.
Beside this, this way of modeling allows us to connect several Territories through their common boundaries and then to consider a large number of connected Territories.\\
Regarding the time scales, in the Toy-Model there are the characteristic time of Population motion, the characteristic time of wealth production and the one of wealth diffusion.

\end{document}